\documentclass[10pt,leqno]{amsart}
\usepackage[a4paper,margin=2cm]{geometry}
\usepackage{MnSymbol}
\newcommand{\EWD}[3]{\left\langle#1:#2:#3\right\rangle}

\def\Card{\mathop{\rm Card}\nolimits}
\def\Char{\mathop{\rm Char}\nolimits}
\def\Coef{\mathop{\rm Coef}\nolimits}

\def\Dom{\mathop{\rm Dom}\nolimits}
\def\Expec{\mathop{\mathcal E}\nolimits}
\def\Fun{\mathop{\rm Fun}\nolimits}

\def\KL{\mathop{\rm Dist_{KL}}\nolimits}

\def\Op{\mathop{\rm Op}\nolimits}

\def\Probab{\mathop{\mathcal P}\nolimits}
\def\Sup{\mathop{\rm Sup}\nolimits}

\def\d{{\mathrm{d}}}
\def\N{{\mathbb{N}}}
\def\R{{\mathbb{R}}}

\begin{document}

\title[Estimating the Missing Mass, ...]{Estimating the Missing Mass, Partition Function or Evidence for a Case of Sampling from a Discrete Set}
\author[B J Braams]{Bastiaan J. Braams}
\date{\today}
\address{Centrum Wiskunde \& Informatica (CWI), Amsterdam, The Netherlands}
\email{b.j.braams@cwi.nl}

\begin{abstract}
We consider the problem of estimating the missing mass, partition function or evidence and its probability distribution in the case that for each sample point in the discrete sample space its (unnormalized) probability mass is revealed.
Estimating the missing mass or partition function (evidence) is a well-studied problem for which, in different contexts, the harmonic mean estimator and the Good-Turing (and related) estimators are available.
For sampling on a discrete set with revealed probability masses these estimators can be Rao-Blackwellized, leading to self-consistent estimators not involving an auxiliary distribution with known total mass.
For the case of sampling from a mixture distribution this offers the perspective of anchoring the estimator at both ends: at the diffuse end (high temperature in statistical physics) via an explicit expression for the total probability mass and at the peaked end (low temperature) via the feature of repeated entries in the sample.

Estimation is model-free, but to provide a probability distribution for the missing mass or partition function a model is needed for the distribution of mass.
We present one such model, identify sufficient reduced statistics, and analyze the model in various ways -- Bayesian, profile likelihood, maximum likelihood and moment matching -- with the objective of eliminating the mathematical (nuisance) parameters for a final expression in terms of the observed data.
The most satisfactory (explicit and transparent) result is obtained by a mixed method that combines Bayesian marginalization or profile likelihood optimization for all but one of the parameters with plain maximum likelihood optimization of the final parameter.

\end{abstract}

\maketitle

\bigskip

\section*{Introduction}

Let $D$ be a finite or countably infinite set on which an unnormalized (but normalizable) probability distribution $p:D\to\R_+$ is defined.
Let $Z=\EWD{\sum i}{i\in D}{p(i)}$.
(Bracket notation following E. W. Dijkstra; see below.)
$Z$ is the partition function in the language of statistical physics or the evidence in the language of Bayesian data analysis.
A sample of points $i\in D$ is drawn independently with replacement from the distribution $p/Z$ and the result is encoded in a vector of integer counts $c:D\to\N$.
Moreover, for the elements $i\in D$ drawn at least once ($1\leq c(i)$) the unnormalized probability mass $p(i)$ is revealed.
Let $S\subseteq D$ be the set of points drawn at least once.
We are concerned here with estimating the missing mass $W=\EWD{\sum i}{i\in D\setminus S}{p(i)}$ and its probability distribution.
When added to the observed mass $V=\EWD{\sum i}{i\in S}{p(i)}$ this provides an estimate and a probability distribution for $Z$.

The problem just described is a variation on the missing mass problem of Good-Turing fame \cite{Good1953,CG1991,GS1995} with the nonstandard feature that the masses $p(i)$ are assumed to be known exactly for every sampled point $i\in S$.
Statisticians may find this to be cheating; in any regular problem in statistics there is an unknown underlying distribution and the task of the analyst is to infer something about that distribution.
The probability mass should not just be given to the analyst, not even one point at a time.

In the context of statistical physics the feature that probability mass is revealed for sampled points is common.
The sampling may be from a Boltzmann distribution in which the energy function and the inverse temperature are known, but the space $D$ is so large that it is necessary to rely on sampling to evaluate any quantity of interest.
Similarly, in the context of Bayesian analysis, $p$ is a prior times a likelihood, both known, and in a high-dimensional setting quantities of interest can only be evaluated by sampling.
Having access to the (unnormalized) probability mass makes it possible to resample the counts and thereby to Rao-Blackwellize any estimator.
I believe that the development of that approach (in Part 1 of this work) is an original contribution to partition function or evidence estimation.

To be clear, in a statistical physics or Bayesian analysis context, if the objective is to evaluate the partition function or the evidence then normally some dedicated sampling strategy will be applied, for example nested sampling \cite{Ski2006,ABB2022}, bridge sampling \cite{Ben1976,MW1996,GSM2017} or any of a variety of other methods \cite{CMG2010,FW2012,LMDLS2023}.
A basic tool within such a strategy is importance sampling, with the sample drawn from an auxiliary probability distribution for which the normalization constant is known.
In the present work the focus is on the analysis phase, given a sample from an unnormalized probability distribution for which the normalization constant is unknown.
We emphasize self-consistent methods for estimating $Z$ and its probability distribution given the data.

Application of these methods may involve a mixture distribution, for example one that is associated with multi-temperature sampling or with Bayesian data assimilation.
Different from and complementary to the approaches that rely on an auxiliary distribution with known normalization constant (normally a very diffuse distribution associated with high temperature in a thermodynamic system), the self-consistent methods rely on narrow peaks in the likelihood (thus on features of the low-temperature components of a thermodynamic mixture distribution).
The present work is a companion and further development of ref.~\cite{Bra2025a,Bra2025b}.

To conclude the introduction I mention some background literature.
For sampling from a discrete set refs \cite{SK2018,Zan2020,KHW2020} are of interest.
The best known examples from statistical physics are the Ising model and the Potts model \cite{RC2004,RK2016,LB2021,Izen2021}.
For Bayesian methods in the physical sciences I only note ref.~\cite{LDT2014}.

\section*{Outline}

After some notational matters we start the technical content with a discussion (in Part 1) of approaches to estimating the missing mass in a model-free manner; i.e., without relying on any likelihood function for the probability distribution $p$.
The celebrated Good-Turing estimator belongs to this class as does the harmonic mean estimator.
We consider adaptations (involving inverse probability weighting or Rao-Blackwellization) that use the additional data in the revealed masses $p(i)$ ($i\in S$).
The key property that makes this additional data valuable (for the case of a discrete sample space) is that some elements may be sampled more than once.
We explain how this may be of special value in connection with sampling from a mixture distribution where (in the context of statistical physics) the elements with repeated draws are associated with low-temperature components while for the highest-temperature component an explicit partition function may be available.

We reflect then on the situation that the Rao-Blackwell procedure relies essentially on the concept of sufficient statistics, but this concept (as it is defined in statistics) makes no sense in the context of sampling from a distribution that is completely specified.
So we are led (in Part 2) to introduce a model by which the original problem from statistical physics or Bayesian analysis, with sampling being used because direct computations are intractable, is replaced by a problem from statistics, where the rules of the game are that we must do inference using sampled data only.
We analyze this problem in multiple ways: Bayesian, profile likelihood, a mixed method and (in appendices) maximum likelihood and moment matching.
The goal for any of these approaches is an explicit representation of the distribution of missing mass, partition function or evidence in terms of the data, after eliminating the nuisance parameters of the mathematical model.
This goal remains somewhat elusive.
The most satisfactory outcome is one in which all but one of the model parameters is eliminated analytically (by marginalization or by profile likelihood maximization) and the final parameter is removed by plain maximum likelihood numerical optimization.

\section*{Notation}

For quantified and similar expressions we use notation associated with E.\ W.\ Dijkstra~\cite{DS2012}.
If $\Op$ is any associative and commutative operation (e.g., sum or product in the present work) then $\EWD{\Op i}{\Dom(i)}{\Fun(i)}$ is the result of applying $\Op$ to the values $\Fun(i)$ as the index or argument $i$ varies over the domain indicated by the boolean expression $\Dom(i)$.
The scope of the dummy variable $i$ is delimited by the angle brackets.
Without any operation, $\EWD{i}{\Dom(i)}{\Fun(i)}$ denotes a function and its domain.
(In these expressions, instead of just $i$ there can be any number of indices or arguments.)
Still following~\cite{DS2012}, if $f$ is a vector, matrix, function, or other structure that involves indices or arguments, and if $\Fun(f)$ is a boolean expression that may be understood pointwise, then $[\Fun(f)]$ means $\Fun(f(.))$ at all points of the domain of $f$.
(Example: $[0\leq f]$ to say that $f$ is everywhere nonnegative.)
Several compatible structures may be involved, as in the expression $[f\leq g]$.

$\Char$ is the characteristic function of a relation: $\Char({\rm False})=0$ and $\Char({\rm True})=1$.
$\Card$ is used for the cardinality of a set.
Subscript $|$ is used to denote the restriction of a function or vector to a subset; so with $p$ defined on $D$ and $S\subseteq D$, $p_{|S}$ is the restriction of $p$ to $S$.
The symbolic expression (not a function) $\Expec()$ denotes an expectation value and similarly $\Probab()$ denotes a probability.

The symbols $D$, $p$, $c$, $S$, $V$, $W$ and $Z$ have already been used and they retain their meaning throughout the paper.
Furthermore, $M$ is always the number of distinct sampled points and $N$ is the sample size including repetitions.
To recap: $D$ is primordial.
In the first part of this paper $p$ is taken as fixed (but unknown a priori), while in the second part it is a random variable drawn from some distribution.
The counts vector $c$ is obtained by sampling, so it is a random variable, and then there are derived quantities
\begin{equation}
\begin{split}
    S&=\{i:i\in D\land1\leq c(i)\}\\
    M&=\EWD{\sum i}{i\in S}{1}\qquad\textrm{($=\Card S$)}\\
    N&=\EWD{\sum i}{i\in S}{c(i)}\\
    V&=\EWD{\sum i}{i\in S}{p(i)}\\
    W&=\EWD{\sum i}{i\in D\setminus S}{p(i)}\\
    Z&=\EWD{\sum i}{i\in D}{p(i)}\qquad\textrm{($=V+W$)}\;.
\end{split}
\end{equation}
All except $Z$ are random variables, and in Part 2 of the paper $Z$ is a random variable too.
I am not able to adhere to the convention of mathematical statistics whereby upper case is used strictly for the random variable while lower case is used for realizations of the random variable or for the variable in a generic context such as a function definition.
Following the conventions of statistical physics I use the same capitalization throughout for any variable.

\section*{Part 1: Estimates}

Part 1 is devoted to methods to estimate $W$ or $Z$ without any model for the distribution $p:D\to\R_+$, especially self-consistent methods that employ sampling from $p/Z$ without involving an auxiliary distribution.
With $D$ a discrete set any point may occur multiple times in the sample, and this points to the Good-Turing estimator and the inverse birthday paradox problem (estimating the number of days in the year).
Different from the classical missing species problem, here we have access to the true (unnormalized) probability masses for observed species,
which invites inverse probability weighting or Rao-Blackwellization.
At the end of Part 1 we discuss application of the self-consistent estimators in connection with continuation methods that may involve an auxiliary distribution with known total mass.

\section*{Sampling framework}

The counts vector $c:D\to\N$ is obtained by sampling independently and with replacement from $p/Z$.
If the sampling is done for a given sample size $N$ then $c$ is drawn from a multinomial distribution
\begin{equation}
    c\sim\EWD{m}{m\in\N^D\land\EWD{\sum i}{i\in D}{m(i)}=N}{\frac{N!}{Z^N}\EWD{\prod i}{i\in D}{\frac{p(i)^{m(i)}}{m(i)!}}}
\end{equation}
We will often assume that the sampling is done instead as a Poisson process with rate parameter $\lambda Z$.
In that case all the $c(i)$ are independent, $N$ is a derived quantity and
\begin{equation}
\begin{split}
    c(i)&\sim\EWD{k}{0\leq k}{e^{-\lambda p(i)}(\lambda p(i))^k/k!}\qquad\textrm{($i\in D$)}\\
    N&\sim\EWD{k}{0\leq k}{e^{-\lambda Z}(\lambda Z)^k/k!}\;.
\end{split}
\end{equation}
Then $\Expec(c(i))=\lambda p(i)$ and $\Expec(N)=\lambda Z$.
For the case of sampling via a Poisson process estimating $\lambda$ given $N$ is essentially the same as estimating $Z$.

Given the unnormalized $p$ together with the normalization constant $Z$ and the sampling protocol (including its parameter $N$ or $\lambda$), $S$ is a random set.
We denote by $\pi:D\to\R_+$ the first-order inclusion probabilities
\begin{equation}
    \pi(i)=\Probab(i\in S)\qquad\textrm{($i\in D$)}\;.
\end{equation}
In the present work, as will be clarified, the precise counts $c$ are generally irrelevant; what matters is the binary property $i\in S$ (for $i\in D$).
From that perspective the sampling strategy looks like sampling without replacement and with inclusion probabilities $\pi$.
However, that leaves higher-order inclusion correlations unspecified.

\section*{Estimators employing inverse probability weighting}

Let $h:D\to\R$ be any quantity that is absolutely summable and let
\begin{equation}
    H=\EWD{\sum i}{i\in D}{h(i)}\;.
\end{equation}
An unbiased estimate of $H$ is obtained by noting that $\Expec(\Char(i\in S))=\Probab(i\in S)$ and therefore
\begin{equation}
\begin{split}
    H&\simeq\EWD{\sum i}{i\in D}{h(i)\Char(i\in S)/\pi(i)}\\
    &=\EWD{\sum i}{i\in S}{\frac{h(i)}{\pi(i)}}\;.
\end{split}
\end{equation}
Note that $[0<p_{|S}]\land[0<p\Rightarrow0<\pi]$, therefore $[0<\pi_{|S}]$.

In survey sampling this estimator for $H$ is often called the (Narain-)Horvitz-Thompson estimator, with reference to \cite{Nar1951,HT1952}.
In fact, those works are devoted to analysis of various sampling strategies (without replacement) and their first-order and higher-order inclusion probabilities.
The basic estimator for $H$ based on weighting by inverse probability of inclusion is taken for granted as just completely classical.
I follow the practice of calling it inverse probability weighting.

In order to obtain a self-consistent estimate of $Z$ we apply inverse probability weighting to the case $h=p$, thus $H=Z$.
\begin{equation}
\begin{split}
    Z&=\EWD{\sum i}{i\in D}{p(i)}\\
    &\simeq\EWD{\sum i}{i\in S}{\frac{p(i)}{\pi(i;Z)}}\;.
\end{split}
\end{equation}
(We emphasize that $\pi$ depends on $Z$, but of course it also depends on $p$ and on the sampling protocol and its parameters.)
Let us spell this out for the two principal sampling strategies considered in this work.

For sampling at fixed $N$, $\Probab(i\in S)=1-(1-p(i)/Z)^N$ ($i\in D$).
Then
\begin{equation}
    Z\simeq\EWD{\sum i}{i\in S}{\frac{p(i)}{1-(1-p(i)/Z)^N}}\;.
\end{equation}
This is exact in expectation (the estimate is unbiased).
Taking the approximation as an identity we have a nonlinear equation for $Z$ of which the solution is an estimator for $Z$.
In the case $M=N$ one must let $Z\to\infty$; there is no finite solution.

If sampling is done as a Poisson process with rate parameter $\lambda Z$ then $\Probab(i\in S)=1-e^{-\lambda p(i)}$ ($i\in D$).
In that case $N$ is an observed quantity and we follow the practice of identifying $\lambda$ with $N/Z$.
Then
\begin{equation}
    Z\simeq\EWD{\sum i}{i\in S}{\frac{p(i)}{1-e^{-Np(i)/Z}}}\;.
\end{equation}
Also in this case, taking the approximation as an identity we have a nonlinear equation for $Z$ of which the solution is an estimator for $Z$.
And here too, if $M=N$ then one must let $Z\to\infty$.

These two estimators for $Z$ (as the solution to a nonlinear equation) both look to me like something that might be found in the 1950s Monte Carlo literature, but I am not able to provide an attribution.
For both estimators there is an implicit assumption that the missing mass is unbounded; they both fail to say $W=0$ in the special case $S=D$.
Indeed, these estimators are completely agnostic about $D$.
The largest difference arises for those points $i\in S$ for which $Np(i)/Z\simeq1$, so points that are expected to be sampled approximately once.
Between the two approaches I favor the fixed-$N$ approach as it avoids the detour over the rate parameter $\lambda$, but I don't have a real conclusive argument.

\section*{Rao-Blackwellized estimators}

Let $f:D\to\R$ be any quantity that is absolutely summable with respect to $p$ and let
\begin{equation}
    \bar f=\frac1{Z}\EWD{\sum i}{i\in D}{p(i)f(i)}\;.
\end{equation}
The standard Monte Carlo estimate for $\bar f$ based on the sample that is represented by the counts vector $c$ is
\begin{equation}
    \bar f\simeq\frac1{N}\EWD{\sum i}{i\in S}{c(i)f(i)}\;.
\end{equation}
For the points $i\in S$ the probability $p(i)$ is revealed.
Therefore, taking $S$, $p_{|S}$ and $N$ as known it is possible to average $c$ over its probability distribution and obtain an estimate for $\bar f$ that has reduced variance.
In the field of Monte Carlo simulation for particle transport this kind of procedure is used in, for example, a track length estimator.
In statistics it is known as the Rao-Blackwell procedure.
For presentations in the context of sampling and Markov Chain Monte Carlo I note especially refs \cite{CR1996,RR2021}.

Given $S$, $p_{|S}$ and $N$ the counts vector $c$ is drawn from a truncated multinomial distribution
\begin{equation}\label{eq:tmc}
    c\sim\EWD{m}{m\in\N^D\land[1\leq m_{|S}]\land[m_{|D\setminus S}=0]\land\EWD{\sum i}{i\in S}{m(i)}=N}{\frac{N!}{F_N(S,p_{|S})}\EWD{\prod i}{i\in S}{\frac{p(i)^{m(i)}}{m(i)!}}}\;.
\end{equation}
Here, $F_N(S,p_{|S})$ is the normalization constant
\begin{equation}
    F_N(S,p_{|S})=\EWD{\sum m_{|S}}{[1\leq m_{|S}]\land\EWD{\sum i}{i\in S}{m(i)}=N}{N!\EWD{\prod i}{i\in S}{\frac{p(i)^{m(i)}}{m(i)!}}}\;.
\end{equation}
The Rao-Blackwellized estimator for $\bar f$ is obtained by replacing $c$ by its expectation value given $(S,p_{|S},N)$.
\begin{equation}
    \bar f\simeq\frac1{N}\EWD{\sum i}{i\in S}{v(i;S)f(i)}\;.
\end{equation}
in which
\begin{equation}
\begin{split}
    v(i;S)&=\Expec(c(i)|S)\\
    &=p(i)\frac\partial{\partial p(i)}\log(F_N(S,p_{|S}))\;.
\end{split}
\end{equation}
Note that $[1\leq v_{|S}]$, $[v_{|D\setminus S}=0]$ and $\EWD{\sum i}{i\in S}{v(i;S)}=N$.
If $M=N$ (all sampled points are singletons) the result is $v_{|S}=\EWD{i}{i\in S}{1}$ and the Rao-Blackwellized estimator is just the original estimator.
The interesting cases for this application of Rao-Blackwell are those in which some points are sampled more than once.

The coefficients $\EWD{i}{i\in S}{v(i;S)}$ may be computed exactly by a dynamic programming recurrence over the variable $i$, although this may not be attractive in applications.
I have not found a closed-form expression for $F_N$ or for $v$.
Resampling $c$ by Markov Chain Monte Carlo could be reasonable in some circumstances.
There are also asymptotic approximations for large $N$; this is discussed briefly in an appendix.

One obtains a Rao-Blackwellized estimator for $Z$ (or rather for $1/Z$) by considering a special case.
Let $f=\pi$, the inclusion probability, and we emphasize that $\pi$ depends on the unknown $Z$.
Then $\bar f\simeq V/Z$ and the estimate is unbiased if the true $Z$ is used.
\begin{equation}
    \frac{V}{Z}\simeq\frac1{N}\EWD{\sum i}{i\in S}{v(i;S)\pi(i;Z)}\;.
\end{equation}
Alternatively we can use $f=\pi/p$ and then $\bar f\simeq M/Z$, again unbiased if the true $Z$ is used.
\begin{equation}
    \frac{M}{Z}\simeq\frac1{N}\EWD{\sum i}{i\in S}{v(i;S)\frac{\pi(i;Z)}{p(i)}}\;.
\end{equation}
Either way, taking the estimate as an identity with $V$ or $M$ observed (and $N$ known or observed), this is a nonlinear equation for $Z$.
I understand these two estimates as Rao-Blackwellized instances of the harmonic mean estimator or reciprocal importance sampling as discussed in a later Section.
Either approach is uninformative about $Z$ if $M=N$; it is only useful if enough probability weight is on points that are expected to be sampled more than once.

\section*{Rao-Blackwellized Good-Turing estimator}

Recall $W=\EWD{\sum i}{i\in D\setminus S}{p(i)}$, the unnormalized missing mass.
The Good-Turing estimator for the normalized missing mass $W/Z$ involves $\Phi_1$, the $k=1$ instance of the more general quantity $\Phi_k$ that is the number of objects sampled exactly $k$ times.
\begin{equation}
    \Phi_k=\EWD{\sum i}{i\in D\land c(i)=k}{1}\qquad\textrm{$(0\leq k)$}\;.
\end{equation}
(Thus $M=\EWD{\sum k}{1\leq k}{\Phi_k}$ and $N=\EWD{\sum k}{1\leq k}{k\Phi_k}$.)
The estimator is $W/Z=\Phi_1/N$.
Consistent with that, $Z$ is estimated as $Z=V/(1-W/Z)$, ${}=VN/(N-\Phi_1)$, and $W=Z(W/Z)$, ${}=V\Phi_1/(N-\Phi_1)$.
These estimates do not involve any assumption about the distribution of masses $\EWD{i}{i\in D}{p(i)}$.

For an exact Rao-Blackwell treatment it is required to evaluate $\Expec(\Phi_1)$ for $c$ drawn from the truncated multinomial distribution in Eq.~\ref{eq:tmc}.
For analytical insight we may employ the Poisson sampling approximation for Rao-Blackwellization.
Therefore, given $S$, $p_{|S}$ and $N$ the parameter $\lambda$ is determined by
\begin{equation}
    N=\EWD{\sum i}{i\in S}{\frac{\lambda p(i)}{1-e^{-\lambda p(i)}}}\;.
\end{equation}
(If $N=M$ then one must let $\lambda\to0$.)
The expected $\Phi_k$ follows:
\begin{equation}
    \Expec(\Phi_k)=\EWD{\sum i}{i\in S}{\frac{(\lambda p(i))^k/k!}{e^{\lambda p(i)}-1}}\;.
\end{equation}
Approximating $\Expec(\Phi_1/N)$ by $\Expec(\Phi_1)/N$ an approximate Rao-Blackwellized normalized missing mass estimator is:
\begin{equation}
    W/Z\simeq\EWD{\sum i}{i\in S}{\frac{\lambda p(i)/N}{e^{\lambda p(i)}-1}}\;.
\end{equation}
It is nicer to express the result directly in terms of the problem data, eliminating the parameter $\lambda$.
To this end, observe that $\lambda$ is defined so that
\begin{equation}
    1=\EWD{\sum i}{i\in S}{\frac{(\lambda p(i)/N)e^{\lambda p(i)}}{e^{\lambda p(i)}-1}}\;.
\end{equation}
Forcing $V/Z+W/Z=1$ one obtains
\begin{equation}
    V/Z=\EWD{\sum i}{i\in S}{\lambda p(i)/N}\;.
\end{equation}
But $V=\EWD{\sum i}{i\in S}{p(i)}$ and so we may identify $\lambda=N/Z$.
Thereby the nonlinear equation for $\lambda$ may be replaced with a nonlinear equation directly for $Z$:
\begin{equation}
    Z=\EWD{\sum i}{i\in S}{\frac{p(i)}{1-e^{-N p(i)/Z}}}\;.
\end{equation}
And then for the missing mass, using $W=Z-V$:
\begin{equation}
    W=\EWD{\sum i}{i\in S}{\frac{p(i)}{e^{N p(i)/Z}-1}}\;.
\end{equation}
The case $N=M$ is singular and then one has to let $Z\to\infty$, $W\to\infty$, $W/Z=1$ and $V/Z=0$.

\section*{Variations on Good-Turing}

The Good-Turing estimator for the missing mass is one in a family of estimators of the general form
\begin{equation}
    W/Z\simeq\EWD{\sum k}{1\leq k}{c_k\Phi_k}
\end{equation}
Estimators of this kind are reviewed and developed in \cite{Pain2025}.
Among the many earlier references I note \cite{GT1956,ET1976,OSW2016} and for recent work see \cite{Schmitz2025,CDH2025}
The Good-Toulmin missing mass estimator has $c_k=-(-1/N)^k$.
Using the same reasoning and approximations as before (including $\Expec(\Phi_k/N)\simeq\Expec(\Phi_k)/N$) this leads to the Rao-Blackwellized Good-Toulmin estimator in the approximation of inverse probability weighting
\begin{equation}
    W/Z\simeq\EWD{\sum i}{i\in S}{\frac{1-e^{-\lambda p(i)/N}}{e^{\lambda p(i)}-1}}
\end{equation}
The Good-Turing problem is actually richer than the missing mass problem; it asks to estimate the mass for any species, seen or unseen.
But in the present context the masses of the seen species are revealed.

\section*{Rao-Blackwellized Harmonic Mean Estimator}

When it is used in isolation the harmonic mean estimator \cite{NR1994} is understood to be a terrible tool to estimate the partition function, and it will be seen that Rao-Blackwellization does nothing to ameliorate its deficiencies.
However, under different guises the harmonic mean estimator has a legitimate role in continuation approaches (such as thermodynamic integration) to estimate a ratio of partition functions for two neighboring distributions.
In such a context the harmonic mean estimator is more often called reciprocal importance sampling \cite{GD1994}, but the mathematical core is the same.
I discuss the harmonic mean estimator here without assuming a particular context.

Suppose that for some function $h:D\to\R$ the unweighted sum is available:
\begin{equation}
    H=\EWD{\sum i}{i\in D}{h(i)}\;.
\end{equation}
For the basic harmonic mean estimator one normally thinks of $h$ as some diffuse positive function; for example with $h(i)$ as a volume element if $D$ is a discretization of a continuous space, or with $[h=1]$ if $D$ is finite.
For application in continuation methods $h$ may be some unnormalized probability density for which the partition function has been determined and for which $h/H$ is close to $p/Z$.

Sampling from the distribution $p/Z$ gives rise to the approximation
\begin{equation}
    H\simeq\frac{Z}{N}\EWD{\sum i}{i\in S}{c(i)h(i)/p(i)}\;.
\end{equation}
(If $Z$ is known then this is the importance weighting estimate for $H$.)
With $H$ and $p_{|S}$ known and $Z$ unknown one obtains the harmonic mean estimator (or reciprocal importance sampling).
\begin{equation}
    Z\simeq NH/\EWD{\sum i}{i\in S}{c(i)h(i)/p(i)}\;.
\end{equation}
If it is computationally feasible then one may Rao-Blackwellize the estimator for $H$ to obtain
\begin{equation}
    H\simeq\frac{Z}{N}\EWD{\sum i}{i\in S}{v(i;S)\frac{h(i)}{p(i)}}\;.
\end{equation}
This is a simple linear equation for $Z$ directly mirroring the classical harmonic mean estimator.
On the other hand, if Rao-Blackwellization is approximated by inverse probability weighting then one obtains
\begin{equation}
    H\simeq\EWD{\sum i}{i\in S}{\frac{h(i)}{\pi(i,Z)}}\;.
\end{equation}
With $H$ known this is a nonlinear equation for $Z$ and the solution is an approximate Rao-Blackwellized harmonic mean estimator.
Both estimators reduce to the original harmonic mean estimator in the limit in which $Np(i)/Z\ll1$ for all $i\in S$ and $[c_{|S}=1]$.
Therefore, as before, this Rao-Blackwellization or inverse probability weighting is of interest in cases in which a significant fraction of the sample is made up of points that are expected to be sampled more than once.

Just to be sure, terms $c(i)h(i)/p(i)$ in the estimator for $H$ for which $p(i)/Z\ll h(i)/H$ cause the harmonic mean estimator to have unbounded variance.
Rao-Blackwellization does nothing to ameliorate that; for those terms it is expected that $c(i)=1$ and in the Rao-Blackwellized estimator for $H$ they also contribute as $h(i)/p(i)$.

\section*{Application to mixed sampling strategies}

So far in Part 1 the nature of the underlying probability distribution $p$ has been left unspecified.
It could be a Boltzmann distribution at one temperature or it could be a profile of Bayesian evidence or something else entirely.
In practical applications, when the aim is to estimate $Z$ then almost always some dedicated sampling strategy will be employed, for example involving multiple temperatures in statistical physics or involving data assimilation in Bayesian analysis.
One may then be analyzing observations from a mixture distribution.
In a way, nothing changes as far as the present analysis is concerned, but nevertheless we lay out here the application of the approaches of Part 1 for the case of analyzing a mixture distribution.

We consider $J$ distributions all defined on the same discrete space $D$ and with unnormalized probability densities $r:D\times\{0,..,J-1\}\to\R_+$.
\begin{equation}
    r=\EWD{i,j}{i\in D\land0\leq j<J}{r(i,j)}\;.
\end{equation}
The normalization constants are
\begin{equation}
    R(j)=\EWD{\sum i}{i\in D}{r(i,j)}\qquad\textrm{($0\leq j<J$)}\;.
\end{equation}
Typically $R(0)$ would be known and the other $R(j)$ would be unknown quantities of interest.

The sampling across the mixture involves weights $w=\EWD{j}{0\leq j<J}{w(j)}$ that are most often determined adaptively in the course of the procedure.
Here we take those weights for granted and known.
The sampling is then from the unnormalized mixture distribution $p:D\to\R_+$,
\begin{equation}
    p(i)=\EWD{\sum j}{0\leq j<J}{r(i,j)w(j)}\qquad\textrm{($i\in D$)}\;.
\end{equation}
The total mass or partition function $Z$ is defined with respect to $p$ as before.
\begin{equation}
\begin{split}
    Z&=\EWD{\sum i}{i\in D}{p(i)}\\
    &=\EWD{\sum i,j}{i\in D\land0\leq j<J}{r(i,j)w(j)}\\
    &=\EWD{\sum j}{0\leq j<J}{R(j)w(j)}\;.
\end{split}
\end{equation}
Sampling is done with replacement and a counts vector $c:D\to\N$ is obtained.
For any point $i$ for which $1\leq c(i)$ the unnormalized mass $p(i)$ and all the $\EWD{j}{0\leq j<J}{r(i,j)}$ become known.
The set $S\subseteq D$ and associated quantities $M$, $N$, $V$ and $W$ are defined as before.

Rao-Blackwellization using $(S,p_{|S},N)$ gives us $v$, concentrated on $S$ and for which $\Expec(c_{|S})=v_{|S}$.
There are also the approximations to $v$ based on inverse probability weighting with probability of inclusion $\pi$ obtained either at fixed $N$ or using the Poisson model.
Thus
\begin{equation}
\begin{split}
    R(j)&\simeq\frac{Z}{N}\EWD{\sum i}{i\in S}{v(i;S)r(i,j)/p(i)}\\
    &\simeq\EWD{\sum i}{i\in S}{\frac{r(i,j)}{\pi(i;Z)}}\;.
\end{split}
\end{equation}
These expressions rely on an approximation for the total mass $Z$, and to estimate that we now have a mix of two strategies.
On the one hand, if $R(0)$ is known or, more generally, if any particular quantity $H=\EWD{\sum i}{i\in D}{h(i)}$ is known, then the mathematics of harmonic mean estimation provides an estimator for $Z$.
On the other hand, if in the weighted mixed sample there is adequate weight on points that are sampled more than once then the self-consistent estimators that use inverse probability weighting are available.
(In the statistical physics context this would be the case if the mixture distribution extends far enough into the low temperature region.)
A weighted mix of the two strategies may rely on the approximation
\begin{equation}
    1\simeq\EWD{\sum i}{i\in S}{\frac{\gamma h(i)/H+(1-\gamma)p(i)/Z}{\pi(i;Z)}}\;.
\end{equation}
A weight $\gamma$ ($0\leq\gamma\leq1$) must be chosen and then by taking the approximation as an identity one obtains a nonlinear equation for $Z$.
How to balance the two terms, or when to put all the weight on a single term, is left to studies more closely tied to specific applications.

\section*{Closing remarks for Part 1}

Three self-consistent estimators have been presented: the estimator directly for $Z$ based on inverse probability weighting and the estimators for $V/Z$ and for $M/Z$ based on Rao-Blackwellization together with the inclusion probabilities.
For each estimator there is also the (perhaps minor) variation in the choice of sampling at fixed total size $N$ versus Poisson sampling.
The numerical implementation of the first estimator, based on inverse probability weighting for $Z$, was explored already in \cite{Bra2025a}.
The approximations relying on Rao-Blackwellization are computationally intensive and in practical cases exact Rao-Blackwellization may not be feasible.
I have to leave it open for now if there are practical cases in which the Rao-Blackwell approach may be preferred over the more straightforward self-consistent inverse probability weighting.

The present work invites applied studies of the self-consistent inverse probability weighting approach together with continuation methods such as multi-temperature sampling in statistical physics.
Classically these multi-temperature approaches are anchored at very high temperature (for which the partition function is available), whereas the present self-consistent approaches make it possible to use an anchor (or a second anchor) at sufficiently low temperature without identifying the zero-temperature limit state or states.

\section*{Intermezzo: Statistical physics and statistics}

The problem discussed here, estimating the missing mass using counts $c(i)$ and revealed masses $p(i)$ ($i\in S$), is motivated by statistical physics, and we run into problems when trying to use the language of mathematical statistics.
In particular, in the standard setting in statistical physics one is working with distributions that are completely known, but the state space is so large that one must rely on sampling to make computations tractable.
In statistics, on the other hand, there is an underlying distribution that is unknown, and the task of the analyst is inference on the basis of the observations.

The language problem comes to the fore when we observe that, given the set $S$ together with the unnormalized probabilities $p_{|S}$ and the total count $N$, the precise observed $c_{|S}$ is irrelevant for the missing mass problem; we can resample $c_{|S}$ given the reduced data $S$, $p_{|S}$ and $N$, as in the previous Rao-Blackwell treatments.
Using the language of statistics we would like to invoke sufficiency of $(S,p_{|S},N)$ for the observed $(S,p_{|S},c_{|S})$.
A precise statement would say that the conditional distribution of the observed $c_{|S}$ given $(S,p_{|S},N)$ does not depend on any unknown parameters governing the distribution of $p$.
In the setting of statistical physics this is a vacuous statement as the quantity $p$ is complete specified.
There are no unknown parameters, and (formally) $S$ and $N$ alone already constitute a sufficient statistic.
In fact, the missing mass $W$ is precisely defined if only $S$ is known, while the total mass $Z$ (the partition function or evidence) is precisely defined without any observations.

In order to reconcile the two perspectives I reformulate the statistical physics problem as a statistics problem.
There is an unknown (parameterized with unknown parameters) probability distribution that governs $p$.
There are observations, based on a known sampling strategy, of selected elements $i\in D$ together with their mass $p(i)$ and count $c(i)$.
These observations define the set $S$ for which $[1\leq c_{|S}]$ and $[c_{|D\setminus S}=0]$.
The task is to infer properties of the underlying distribution using only the sampled data.
There is no longer any reference to the state space being too large for explicit computation; instead, the analyst simply does not have access to $p(i)$ outside the sample.
Having replaced intractability (in statistical physics) by ignorance (in statistics) we should no longer assume that the state space is extremely large, as was done in Part 1.
$D$ can be of any size, and even the extreme case $S=D$ (there is no missing mass) must receive a natural treatment.

\section*{Part 2: Distributions}

In Part 1 the unnormalized probability distribution $p$ was a fixed unknown quantity and we considered estimators for the missing mass $W$ or the total mass $Z$.
Bootstrap resampling can provide a distribution of the estimated $W$ or $Z$, but it does not provide a probability distribution for the underlying true value.

In Part 2 we treat $p$ as a random variable drawn from some parametric distribution, for which we make one particular choice.
Given the parametric distribution for $p$ (with unknown parameters) one can obtain a probability distribution for $W$ from the observed data.
Several strategies are available: fully Bayesian analysis, maximum likelihood estimation, profile likelihood estimation, different moment matching approaches and a multitude of mixed methods.

It is of interest to study how the obtained distribution for $W$ depends on the model assumptions and on the mode of analysis.
Ideally model parameters (in the assumed probability distribution for $p$) can be eliminated and the results expressed in terms of the observed quantities $(S,p_{|S},N)$ in some transparent manner.
We look in all directions and find Bayesian estimation and profile likelihood estimation to be the most satisfactory starting points.
A particular mix of either Bayesian or profile likelihood estimation (for all but one parameter) together with maximum likelihood (for the remaining parameter) provides the most natural interpretable results.

\section*{Model specification}

The statistical model for $p$ will be stated first, followed by some justification.

The sample space is a finite or countably infinite set $D$ on which there is a known normalized distribution $x:D\to\R_+$.
(So $\EWD{\sum i}{i\in D}{x(i)}=1$.)
There are unknown positive real parameters $\alpha$, $b$ and $\lambda$.
For each $i\in D$ a nonnegative value $p(i)$ is drawn from the Gamma distribution with shape parameter $a(i)=\alpha x(i)$ and rate parameter (inverse scale parameter) $b$, thus mean value $a(i)/b$,
\begin{equation}
    p(i)\sim\EWD{t}{0\leq t}{t^{a(i)-1}e^{-bt}}\frac{b^{a(i)}}{\Gamma(a(i))}\;.
\end{equation}
(If $a(i)=0$ then $p(i)=0$ almost surely.)
Given $p(i)$ an integer count $c(i)$ is drawn from a Poisson distribution with rate parameter (mean value) $\lambda p(i)$,
\begin{equation}
    c(i)\sim\EWD{k}{0\leq k}{(\lambda p(i))^k/k!}/e^{\lambda p(i)}\;.
\end{equation}
(If $p(i)=0$ then $c(i)=0$ almost surely.)
After the draw, if $c(i)=0$ then only $c(i)$ is revealed.
If $1\leq c(i)$ then $p(i)$ and $c(i)$ are revealed.
Let $S\subseteq D$ be the set $S=\{i:1\leq c(i)\}$ and then let $V=\EWD{\sum i}{i\in S}{p(i)}$, $W=\EWD{\sum i}{i\in D\setminus S}{p(i)}$ and $Z=\EWD{\sum i}{i\in D}{p(i)}$ ($=V+W$).
$V$ is the observed mass, $W$ is the hidden mass and $Z$ is the total mass.
Two additional important observed quantities are $M=\EWD{\sum i}{i\in S}{1}$, ${}=\Card S$, and $N=\EWD{\sum i}{i\in S}{c(i)}$, the sample size including repetitions.
The task for the analyst is to estimate $W$ (or $Z$) and its probability distribution given only the revealed data.

A rationale for this model follows.
At first we might like to say that each $p(i)$ is drawn from some totally uninformative improper prior $\EWD{t}{0\leq t}{1/t}$, but it is not possible to calculate with that prior.
The Gamma distribution is chosen as a related proper prior in which the shape parameter ensures integrability for $p(i)\to 0$ and the rate parameter ensures integrability for $p(i)\to\infty$.
To cater to the situation that $D$ may be an infinite set and also the situation that $D$ may have some structure, we let the shape parameter depend on the point $i\in D$ through the coefficients $x=\EWD{i}{i\in D}{x(i)}$.
Then the parameters $\alpha$ and $b$ together with $x$ control the distribution of $p$ while the parameter $\lambda$ controls the Poisson sampling of $c$.

The model is introduced for the purpose of mathematical analysis and there is no profound meaning to be assigned to the unknown coefficients $(\alpha,b,\lambda)$.
They are nuisance parameters, to be eliminated in the course of the analysis (by marginalization, optimization or moment matching), so that the results are finally expressed in terms of the observations.
Then we would like to see that those end results are interpretable and robust against different model choices.

\section*{Preliminary likelihood analysis of the model}

The joint distribution of the vectors $p$ and $c$ (the likelihood function) is a direct product of Poisson-Gamma distributions.
\begin{equation}
    L0(p,c;\alpha,b,\lambda)=\EWD{\prod i}{i\in D}{\frac{b^{a(i)}}{\Gamma(a(i))}p(i)^{a(i)-1}e^{-(b+\lambda)p(i)}(\lambda p(i))^{c(i)}/c(i)!}\;.
\end{equation}
(If $a(i)=0$ then the associated factor in $L0$ is a unit mass at $p(i)=0\land c(i)=0$, i.e., a factor $\delta_0(p(i))\Char(c(i)=0)$.)
Here and throughout we let $a(i)=\alpha x(i)$.

In a trivial way the likelihood function $L0(p,c;\alpha,b,\lambda)$ can be written as a likelihood involving also the set $S$, and then it may be split into factors involving $i\in S$ and factors involving $i\in D\setminus S$.
\begin{equation}
\begin{split}
    L1(S,p,c;\alpha,b,\lambda)&=\Char(\EWD{\forall i}{i\in D}{i\in S\Leftrightarrow1\leq c(i)})L0(p,c;\alpha,b,\lambda)\\
    &=\EWD{\prod i}{i\in S}{\Char(1\leq c(i))\frac{b^{a(i)}}{\Gamma(a(i))}p(i)^{a(i)-1}e^{-(b+\lambda)p(i)}(\lambda p(i))^{c(i)}/c(i)!}\\
    &\quad\times\EWD{\prod i}{i\in D\setminus S}{\Char(c(i)=0)\frac{b^{a(i)}}{\Gamma(a(i))}p(i)^{a(i)-1}e^{-(b+\lambda)p(i)}}\;.
\end{split}
\end{equation}
We recall the derived quantities $M$, $N$, $V$, $W$ and $Z$ introduced earlier independent of a model for $p$.
In connection with the present model several additional derived quantities are introduced.
\begin{equation}
\begin{split}
    T&=\EWD{\sum i}{i\in S}{x(i)\log x(i)}\\
    U&=\EWD{\sum i}{i\in S}{x(i)\log p(i)}\\
    X&=\EWD{\sum i}{i\in S}{x(i)}\\
    Y&=\EWD{\sum i}{i\in D\setminus S}{x(i)}\qquad\textrm{(so $X+Y=1$)}\;.
\end{split}
\end{equation}
With use of these quantities the likelihood function $L1$ becomes
\begin{equation}
\begin{split}
    L1(S,p,c;\alpha,b,\lambda)&=\EWD{\prod i}{i\in S}{\Char(1\leq c(i))\frac{p(i)^{c(i)-1}/c(i)!}{\Gamma(a(i))}}b^{\alpha X}\lambda^Ne^{\alpha U-(b+\lambda)V}\\
    &\quad\times\EWD{\prod i}{i\in D\setminus S}{\Char(c(i)=0)\frac{p(i)^{a(i)-1}}{\Gamma(a(i))}}b^{\alpha Y}e^{-(b+\lambda)W}\;.
\end{split}
\end{equation}
($T$ is not yet used, but it will come up in some approximations for the case $\alpha\to\infty$.)

We now marginalize over the unobserved quantities $p(i)$ ($i\in D\setminus S$), but retaining the total unobserved mass $W$ as a variable.
A Beta function identity is used.
For a vector $a:D\to\R_+$ and with $\alpha=\EWD{\sum i}{i\in D}{a(i)}$
\begin{equation}
    \EWD{\sumint p}{[0\leq p]\land\EWD{\sum i}{i\in D}{p(i)}=Z}{\EWD{\prod i}{i\in D}{\frac{p(i)^{a(i)-1}}{\Gamma(a(i))}}\,\d p}=\frac{Z^{\alpha-1}}{\Gamma(\alpha)}\;.
\end{equation}
Eliminating the unobserved $p(i)$ from $L1$ by marginalization subject to $\EWD{\sum i}{i\in D\setminus S}{p(i)}=W$ leads to
\begin{equation}
    L2(S,p_{|S},c_{|S},W;\alpha,b,\lambda)=\EWD{\prod i}{i\in S}{\frac{p(i)^{c(i)-1}/c(i)!}{\Gamma(a(i))}}\frac{W^{\alpha Y-1}}{\Gamma(\alpha Y)}b^\alpha\lambda^Ne^{\alpha U-(b+\lambda)(V+W)}\;.
\end{equation}
(If $Y=0$ then the factor $W^{\alpha Y-1}/\Gamma(\alpha Y)$ becomes $\delta_0(W)$; a unit mass at $W=0$.)
Eliminating the unobserved $p(i)$ from $L1$ by marginalization without any constraint leads to
\begin{equation}
    L3(S,p_{|S},c_{|S};\alpha,b,\lambda)=\EWD{\prod i}{i\in S}{\frac{p(i)^{c(i)-1}/c(i)!}{\Gamma(a(i))}}b^\alpha\lambda^Ne^{\alpha U-(b+\lambda)V}(b+\lambda)^{-\alpha Y}\;.
\end{equation}
It is understood that $[1\leq c_{|S}]$.
Indeed, $\EWD{\sumint W}{0\leq W}{L2(S,p_{|S},c_{|S},W;\alpha,b,\lambda)\,\d W}=L3(S,p_{|S},c_{|S};\alpha,b,\lambda)$.

Note that the factor $\EWD{\prod i}{i\in S}{p(i)^{c(i)-1}/c(i)!}$ depends only on the observations $(p_{|S},c_{|S})$.
This factor will cancel out of all normalized probabilities for $(\alpha,b,\lambda)$ or for $W$.
Within the present model the observed $(S,U,V,N)$ are sufficient for the observed $(S,p_{|S},c_{|S})$.
Due to the complicated expression $\EWD{\prod i}{i\in S}{\Gamma(a(i))}$ it is not possible in general to reduce $S$ to a small number of characteristic quantities.
(I won't pursue it, but if $\alpha\to\infty$ and $\Gamma(a(i))$ is replaced by its leading asymptotic approximation then in the tuple of sufficient observations $S$ can be replaced by only $(M,T,X)$.)

\section*{Special cases}

The vector $x_{|S}/X$ is a properly normalized probability distribution on $S$ as is the vector $p_{|S}/V$.
The scaled Kullback-Leibler divergence $\Delta_S=X\KL(x_{|S}/X\|p_{|S}/V)$ will be of interest.
Explicitly
\begin{equation}
\begin{split}
    \Delta_S&=\EWD{\sum i}{i\in S}{x(i)(\log(x(i)/X)-\log(p(i)/V))}\\
    &=T-X\log X-U+X\log V\;.
\end{split}
\end{equation}
In dimension $\Card S+1$ there are the normalized probability distributions $x_{|S}\oplus(Y)$ and $Z^{-1}(p_{|S}\oplus(W))$.
We denote by $\Delta$ their Kullback-Leibler divergence:
\begin{equation}
\begin{split}
    \Delta&=T+Y\log Y-U-Y\log W+\log(V+W)\\
    &=\Delta_S+X\log X+Y\log Y-X\log(V/Z)-Y\log(W/Z)\;.
\end{split}
\end{equation}
The inequalities $0\leq\Delta_S\leq\Delta$ hold.
The special case of equality $\Delta_S=0$ occurs if $p_{|S}$ is proportional to $x_{|S}$; that is, if for some positive $r$, $p_{|S}=x_{|S}r$.
In that case $U=T+X\log r$ and $V=Xr$.
The special case of equality $\Delta=\Delta_S$ occurs if for some positive $r$, $V=Xr$ and $W=Yr$ (and then $Z=r$).
We anticipate that the case $\Delta_S=0$ is going to be somewhat special in the analysis.
If all the observed $p(i)$ are proportional to $x(i)$ with the same constant of proportionality $r$ then a statistical procedure may yield the result that $W=Yr$ almost surely; i.e., the probability distribution for the missing mass becomes a delta function at $Yr$.

The case $M=N$ was special for model-free estimation of the missing mass, but it is not a special case (unless also $\Delta_S=0$) under the given model for $(p,c)$.

\section*{Outlook}

Given the likelihood function there are several canonical ways to proceed.
In the fully Bayesian approach one takes the likelihood $L2(..,W;\alpha,b,\lambda)$ (with $..$ denoting the observations) as the starting point, prescribes a prior distribution on $(\alpha,b,\lambda)$ and marginalizes over these parameters.
The result is an unnormalized distribution for $W$ given the observed data.
The profile likelihood approach also starts from $L2$.
The parameters $(\alpha,b,\lambda)$ are optimized as a function of $W$ and then $L2$ with the optimized parameters is interpreted as an unnormalized distribution for $W$.
On the other hand, the parameters $(\alpha,b,\lambda)$ may be determined from the likelihood $L3(..;\alpha,b,\lambda)$ (not involving $W$) by maximum likelihood estimation or by some form of moment matching.
These parameters are then inserted in $L2$ to obtain an unnormalized distribution for $W$.
For the present problem I found Bayesian estimation and profile likelihood estimation the most appealing among the canonical approaches, but I also pursue a mixed approach using Bayesian estimation or profile likelihood for the parameters $b$ and $\lambda$ and then eliminating $\alpha$ by plain maximum likelihood.
The two other canonical approaches, maximum likelihood (on all parameters) and moment matching, are explored in two Appendices.

\section*{Bayesian estimation of the missing mass}

The starting point for Bayesian estimation of $W$ and its distribution is the likelihood $L2$ together with a prior on the parameters $(\alpha,b,\lambda)$.
The strategy is to marginalize $L2$ over $(\alpha,b,\lambda)$ to give (see below) $L6(S,p_{|S},c_{|S},W)$ and to marginalize $L3$ over $(\alpha,b,\lambda)$ (or integrate $L6$ over $W$) to give $L7(S,p_{|S},c_{|S})$.
Then $L6(..,W)/L7(..)$ is the properly normalized Bayesian probability distribution for $W$ given the observations $(S,p_{|S},c_{|S})$ (and with $[c_{|D\setminus S}=0]$ understood).
The distribution will depend on those observations only through the reduced quantities $(S,U,V,N)$.

The parameters $\alpha$, $b$ and $\lambda$ are nuisance parameters with independent interpretations.
They are all nonnegative real and for the prior I don't know anything better than a non-informative $(\alpha b\lambda)^{-1}$.
With that prior the parameters $b$ and $\lambda$ can be integrated out analytically.

Marginalizing $L2$ and $L3$ over $b$ and $\lambda$ with their prior $(b\lambda)^{-1}$ results in:
\begin{equation}
\begin{split}
    L4(S,p_{|S},c_{|S},W;\alpha)&=\EWD{\prod i}{i\in S}{\frac{p(i)^{c(i)-1}/c(i)!}{\Gamma(a(i))}}\frac{W^{\alpha Y-1}}{\Gamma(\alpha Y)}e^{\alpha U}\Gamma(\alpha)\Gamma(N)(V+W)^{-\alpha-N}\\
    L5(S,p_{|S},c_{|S};\alpha)&=\EWD{\prod i}{i\in S}{\frac{p(i)^{c(i)-1}/c(i)!}{\Gamma(a(i))}}e^{\alpha U}\Gamma(\alpha)\Gamma(N)V^{-\alpha X-N}\Gamma(\alpha X+N)/\Gamma(\alpha+N)\;.
\end{split}
\end{equation}
(If $Y=0$ then the factor $W^{\alpha Y-1}/\Gamma(\alpha Y)$ in $L4$ is understood as $\delta_0(W)$.)

One may check that $\EWD{\sumint W}{0\leq W}{L4(S,p_{|S},c_{|S},W;\alpha)\,\d W}=L5(S,p_{|S},c_{|S};\alpha)$ by using the Beta function identity
\begin{equation}
    \EWD{\sumint W}{0\leq W}{(V+W)^{-\alpha-N}W^{\alpha Y-1}\,\d W}=V^{-\alpha X-N}\frac{\Gamma(\alpha X+N)\Gamma(\alpha Y)}{\Gamma(\alpha+N)}\;.
\end{equation}

Integration of $L4$ and of $L5$ over $\alpha$ (with prior $1/\alpha$) can only be done numerically, but we must study the integrand to convince ourselves that it can be done at all.
That analysis will be found in an Appendix.
It is shown there that the integrand is logarithmically concave and with behavior as $\alpha\to0$ and as $\alpha\to\infty$ that makes it integrable except when $\Delta_S=0$, and the case $\Delta_S=0$ can be treated as a limit.
Then we let
\begin{equation}
\begin{split}
    L6(S,p_{|S},c_{|S},W)&=\EWD{\sumint\alpha}{0\leq\alpha}{L4(S,p_{|S},c_{|S},W;\alpha)\alpha^{-1}\,\d\alpha}\\
    L7(S,p_{|S},c_{|S})&=\EWD{\sumint\alpha}{0\leq\alpha}{L5(S,p_{|S},c_{|S};\alpha)\alpha^{-1}\,\d\alpha}\;.
\end{split}
\end{equation}
The ratio $L6(..,W)/L7(..)$ is the normalized probability distribution for $W$ in the given model for observations $(S,p_{|S},c_{|S})$, but in reality depending only on $(S,U,V,N)$.
In the singular case $\Delta_S=0$ it becomes a delta function $\delta_{Yr}(W)$ in which $r$ is the coefficient of proportionality in $p_{|S}=x_{|S}r$.

The one-dimensional integrals over $\alpha$ in the expressions for $L6$ and $L7$ (with parametric dependence on $W$) are simple enough from the perspective of numerical analysis, but the final expression for $L6(..,W)/L7(..)$ is not very transparent as regards the dependence of $W$ and its distribution on the observed parameters.
The same feature appears in the profile likelihood analysis in the next Section, and we defer further discussion until after that analysis.

\section*{Profile Likelihood Estimation of the missing mass}

The starting point for profile likelihood estimation of $W$ and its distribution is again the likelihood $L2$.
Different from Bayesian estimation no prior is involved for the parameters $(\alpha,b,\lambda)$.
The strategy is to maximize $L2$ over $(\alpha,b,\lambda)$ to give (see below) $L10(S,p_{|S},c_{|S},W)$, which is an unnormalized probability density for $W$ given the observations $(S,p_{|S},c_{|S})$.

Maximization over $b$ and $\lambda$ can be done analytically.
One finds $N=\lambda\alpha/b$ and then
\begin{equation}
\begin{split}
    b(\alpha,W)&=\alpha/(V+W)\\
    \lambda(\alpha,W)&=N/(V+W)\;.
\end{split}
\end{equation}
That leads to the profile likelihood function for $W$ with $\alpha$ as the only parameter.
\begin{equation}
\begin{split}
    L8(S,p_{|S},c_{|S},W;\alpha)&=L2(S,p_{|S},c_{|S},W;\alpha,b(W,\alpha),\lambda(W,\alpha))\\
    &=\EWD{\prod i}{i\in S}{\frac{p(i)^{c(i)-1}/c(i)!}{\Gamma(a(i))}}\frac{W^{\alpha Y-1}}{\Gamma(\alpha Y)}e^{\alpha U}\alpha^\alpha N^N e^{-\alpha-N}(V+W)^{-\alpha-N}\;.
\end{split}
\end{equation}
For later use we record the integral of $L8$ over $W$
\begin{equation}
\begin{split}
    L9(S,p_{|S},c_{|S};\alpha)&=\EWD{\sumint W}{0\leq W}{L8(S,p_{|S},c_{|S},W;\alpha)\,\d W}\\
    &=\EWD{\prod i}{i\in S}{\frac{p(i)^{c(i)-1}/c(i)!}{\Gamma(a(i))}}e^{\alpha U}\alpha^\alpha N^N e^{-\alpha-N}V^{-\alpha X-N}\Gamma(\alpha X+N)/\Gamma(\alpha+N)\;.
\end{split}
\end{equation}
(Note the similarity, especially at large $\alpha$, of $L8$ and $L9$ with the earlier $L4$ and $L5$ of the Bayesian approach.)

Maximization of $L8$ over $\alpha$ (for the given observations and parametrically with respect to $W$) can only be done numerically, but we need to study the function to make sure that it is a well-defined problem.
This analysis is deferred to an Appendix.
It is seen there that $\EWD{\alpha}{0\leq\alpha}{L8(..,W;\alpha)}$ is logarithmically concave and that it has a unique maximum at finite argument except in the singular case $\Delta=0$ when one must let $\alpha\to\infty$.
Therefore we may define
\begin{equation}
    L10(S,p_{|S},c_{|S},W)=\EWD{\Sup\alpha}{0<\alpha}{L8(S,p_{|S},c_{|S},W;\alpha)}\;.
\end{equation}
As is the case for Bayesian estimation, the resulting expression for the profile likelihood of $W$ given the observations is simple enough from the perspective of numerical analysis, but it is not very transparent.
This leads us to a mixed approach that is presented in the next Section.

\section*{A mixed treatment of the missing mass}

In the given model Bayesian estimation and profile likelihood estimation can be carried through analytically to eliminate the parameters $b$ and $\lambda$, but marginalization or profile maximization over $\alpha$ can only be done numerically.
However, the likelihood functions $L4(S,p_{|S},c_{|S},W;\alpha)$ and $L8(S,p_{|S},c_{|S},W;\alpha)$ show a quite simple dependence on $W$ at fixed $\alpha$ and the normalized likelihood functions $L4(..,W;\alpha)/L5(..;\alpha)$ and $L8(..,W;\alpha)/L9(..;\alpha)$ are identical.
In both cases one finds a Beta prime distribution for $W$ or for the scaled variable $W/V$.
\begin{equation}\label{eq:mixa}
\begin{split}
    W&\sim\frac{V^{\alpha X+N}\Gamma(\alpha+N)}{\Gamma(\alpha Y)\Gamma(\alpha X+N)}\EWD{t}{0\leq t}{t^{\alpha Y-1}(V+t)^{-\alpha-N}}\\
    W/V&\sim\frac{\Gamma(\alpha+N)}{\Gamma(\alpha Y)\Gamma(\alpha X+N)}\EWD{t}{0\leq t}{t^{\alpha Y-1}(1+t)^{-\alpha-N}}\;.
\end{split}
\end{equation}
This can be massaged into an ordinary Beta distribution for the normalized missing mass $W/Z$.
\begin{equation}\label{eq:mixb}
    W/Z\sim\frac{\Gamma(\alpha+N)}{\Gamma(\alpha Y)\Gamma(\alpha X+N)}\EWD{t}{0\leq t\leq1}{t^{\alpha Y-1}(1-t)^{\alpha X+N-1}}\;.
\end{equation}
The associated expectation values are
\begin{equation}
\begin{split}
    \Expec(W/V)&=\frac{\alpha Y}{\alpha X+N-1}\\
    \Expec(W/Z)&=\frac{\alpha Y}{\alpha+N}\;.
\end{split}
\end{equation}
From a mathematical point of view, and with an interest in interpretable results, it seems attractive to combine the Bayesian or profile likelihood treatment of the parameters $b$ and $\lambda$ with a plain maximum likelihood estimation of $\alpha$ based upon either $L5(..;\alpha)$ or $L9(..;\alpha)$.
The maximum likelihood estimation of $\alpha$ is still a numerical task, but at least the expression for the distribution of $W$ is then completely explicit.

Perhaps it is also possible to develop a story that never involves the parameters $b$ and $\lambda$ but that leads directly to the likelihood function or probability density of $W$ or of $W/Z$ depending only on the single parameter $\alpha$ (and on the a priori probabilities $x$ and the observations).
I leave that open for now.

\section*{Other approaches to estimating the missing mass}

The previous Sections described the application of Bayesian estimation and profile likelihood estimation, and of a mixed approach in which the parameter $\alpha$ (only) is treated by regular maximum likelihood.
Plain maximum likelihood estimation is another canonical approach, and there are also various forms of moment matching to determine the three parameters.
The analysis for maximum likelihood estimation and instances of moment matching is presented in the Appendices.
For this particular problem I find both those approaches less satisfactory (in the analysis, not for any reasons of principle) than either Bayesian analysis or profile likelihood.
The analysis is there so that readers can make their own assessments.

\section*{Closing remarks on Part 2}

The desired outcome of the analysis in Part 2, employing a model distribution for $p$, is a probability distribution for $W$ (or for $Z$ or for $W/Z$) in terms of the (sufficient) observed quantities $(S,N,U,V)$ after eliminating the model parameters.
The dependence of the distribution of $W$ on the observations should be interpretable and natural, and the end result should not depend too much on the model for $p$, which may have been selected primarily for mathematical convenience.
With that objective in mind the most satisfactory result is that of the mixed treatment, Eqs.~\ref{eq:mixa} and \ref{eq:mixb}.
However, these expressions still involve the parameter $\alpha$ that needs to be determined by numerical optimization, leaving the result not entirely transparent in terms of the observations.

I also have to leave it unresolved for now how closely this outcome is tied to the specific model for $p$.
The likelihood function $L2$ used in the present work is a direct product of Gamma distributions, and a discrete base measure $x:D\to\R_+$ is involved.
That base measure could also be used in a Pitman-Yor process \cite{PY1997} to define a likelihood function.
The Pitman-Yor process has been invoked in connection with species sampling or missing mass problems, e.g.~in \cite{LPR2020,BFN2022}.
There are other approaches to the missing mass problem that are aimed (as here in Part 2) at the distribution of mass, including refs~\cite{OS2015,FNT2016,HNWSW2025,Mour2025}.
I have not tried to connect those approaches to the present problem setting that includes revealed probability mass.

\section*{Conclusions}

In Part 1, estimation, three self-consistent estimators were presented: an estimator directly for $Z$ based on inverse probability weighting and estimators for $V/Z$ and for $M/Z$ based on Rao-Blackwellization together with the inclusion probabilities.
We discussed how these estimators may be of special value in connection with sampling from a mixture distribution where (in the context of statistical physics) the elements with repeated draws are associated with low-temperature components while for the highest-temperature component an explicit partition function may be available.
This obviously invites further study in connection with specific applications.

In Part 2, distributions, the most satisfactory (explicit and transparent) result was obtained by a mixed method that combines Bayesian marginalization or profile likelihood optimization for all but one of the parameters with plain maximum likelihood optimization of the final parameter.
Here the main open question (I think) is to what extent the results are robust against different model assumptions; for example, what kind of distribution of the missing mass would result from a Pitman-Yor model, in the present context where probability mass is revealed for any sample point.

\section*{Acknowledgments}

I thank Wouter Koolen, Peter Grünwald and Charles Karney for valuable discussions.

\section*{Appendix: Remarks on notation}

In Dijkstra's bracket notation the scope of the dummy variables is delimited by the brackets, and so a dummy variable inside angular brackets is distinct from any variable with the same name outside those brackets.
Therefore it is fine to express a normalized probability at the point $i\in D$ as
\begin{equation}
    p(i)/\EWD{\sum i}{i\in D}{p(i)}\;.
\end{equation}
Likewise, for a variable $Z$ drawn from a Gamma distribution one may write
\begin{equation}
    Z\sim\frac{b^\alpha}{\Gamma(\alpha)}\EWD{Z}{0\leq Z}{Z^{\alpha-1}e^{-bZ}}\;.
\end{equation}
If one has internalized the scoping rules then one may appreciate this standardization on $i$ for an element of the state space or on $Z$ for the partition function.

In the square bracket notation for the ``everywhere'' operator the expression inside the brackets must be understood pointwise.
So for compatible numerical structures $f$ and $g$, the value of $f\leq g$ is a boolean structure.
If one wants to have an order operator on structures then a different notation must be used; say $f\preceq g$ or $f\sqsubseteq g$ with a defined operator $\preceq$ or $\sqsubseteq$.
This works fine for me except for equality, where $f=g$ really must mean structure equality; it is not a pointwise expression of which the value is a boolean structure.
For boolean operands Dijkstra resolves this by using $\equiv$ or $\Leftrightarrow$ as the pointwise equivalence operator and using $=$ like everyone else for equality in general.
For numerical structures there is not such an established treatment.
(One could let $\equiv$ mean equality as a pointwise expression for any base type, but I am not going there.)
In conclusion, it is $f=g$ for complete equality between compatible structures of any base type.

\section*{Appendix: Alternative presentations of the model}

Some inspiration may be had from alternative views of the model, all leading to the same observed variables and the same likelihood $L0$, or the likelihood $L2$ after marginalizing away the irrelevant variables.

We described the masses $\EWD{i}{i\in D}{p(i)}$ as being drawn independently from Gamma distributions that share the rate parameter $b$.
In an equivalent description, first $Z$ is drawn from a Gamma distribution and then $p$ is drawn from a Dirichlet (or Beta) distribution with constrained total mass.
\begin{equation}
\begin{split}
    Z&\sim\frac{b^\alpha}{\Gamma(\alpha)}\EWD{t}{0\leq t}{t^{\alpha-1}e^{-bt}}\\
    p&\sim\frac{\Gamma(\alpha)}{Z^{\alpha-1}}\EWD{u}{u\in\R_+^D\land\EWD{\sum i}{i\in D}{u(i)}=Z}{\EWD{\prod i}{i\in D}{\frac{u(i)^{a(i)-1}}{\Gamma(a(i))}}}\;.
\end{split}
\end{equation}
(As always, $a(i)=\alpha x(i)$.)
Similarly, given $p$ the counts $\EWD{i}{i\in D}{c(i)}$ may be obtained by first drawing $N$ from a Poisson distribution and then drawing $c$ from a multinomial distribution at fixed total count $N$.
\begin{equation}
\begin{split}
    N&\sim\EWD{k}{0\leq k}{(\lambda Z)^k/k!}/e^{\lambda Z}\\
    c&\sim\EWD{m}{m\in\N^D\land\EWD{\sum i}{i\in D}{m(i)}=N}{\frac{N!}{Z^N}\EWD{\prod i}{i\in D}{\frac{p(i)^{m(i)}}{m(i)!}}}\;.
\end{split}
\end{equation}

One can also invert the order of drawing $p$ and drawing $c$.
In that case, first $c$ is drawn from its marginal of $L0$, which is a negative binomial distribution for each $c(i)$ independently.
\begin{equation}
\begin{split}
    c(i)&\sim\frac{b^{a(i)}}{\Gamma(a(i))}\EWD{k}{0\leq k}{\EWD{\sumint t}{0\leq t}{t^{a(i)+k-1}e^{-(b+\lambda)t}\,\d t}\lambda^k/k!}\\
    &=\frac{b^{a(i)}}{\Gamma(a(i))}\EWD{k}{0\leq k}{\Gamma(a(i)+k)(b+\lambda)^{-a(i)-k}\lambda^k/k!}\\
    &=\EWD{k}{0\leq k}{\frac{\Gamma(a(i)+k)}{\Gamma(a(i))k!}\left(\frac{b}{b+\lambda}\right)^{a(i)}\left(\frac\lambda{b+\lambda}\right)^k}\qquad\textrm{($i\in D)$}\;.
\end{split}
\end{equation}
The associated probability density for $N$ (still before drawing $p$) is
\begin{equation}
    N\sim\EWD{k}{0\leq k}{\frac{\Gamma(\alpha+k)}{\Gamma(\alpha)k!}\left(\frac{b}{b+\lambda}\right)^\alpha\left(\frac\lambda{b+\lambda}\right)^k}\;.
\end{equation}
The distribution of $\EWD{i}{i\in D}{c(i)}$ given $N$ (but before drawing $p$) has a particularly simple form, being a multinomial distribution that is independent of the parameters $(\alpha,b,\lambda)$.
\begin{equation}
    c\sim\EWD{m}{m\in\R_+^D\land\EWD{\sum i}{i\in D}{m(i)}=N}{N!\EWD{\prod i}{i\in D}{\frac{x(i)^{m(i)}}{m(i)!}}}\;.
\end{equation}

When drawing $c$ before $p$ then each $p(i)$ is drawn independently from a Gamma distribution.
\begin{equation}
    p(i)\sim\frac{(b+\lambda)^{a(i)+c(i)}}{\Gamma(a(i)+c(i))}\EWD{t}{0\leq t}{t^{a(i)+c(i)-1}e^{-(b+\lambda)t}}\;.
\end{equation}
The associated probability densities for $V$, $W$ and $Z$ are
\begin{equation}
\begin{split}
    V&\sim\frac{(b+\lambda)^{\alpha X+N}}{\Gamma(\alpha X+N)}\EWD{t}{0\leq t}{t^{\alpha X+N-1}e^{-(b+\lambda)t}}\\
    W&\sim\frac{(b+\lambda)^{\alpha Y}}{\Gamma(\alpha Y)}\EWD{t}{0\leq t}{t^{\alpha Y-1}e^{-(b+\lambda)t}}\\
    Z&\sim\frac{(b+\lambda)^{\alpha+N}}{\Gamma(\alpha+N)}\EWD{t}{0\leq t}{t^{\alpha+N-1}e^{-(b+\lambda)t}}\;.
\end{split}
\end{equation}
If $V$ is drawn first (after drawing $c$) then the density for the visible mass $p_{|S}$ becomes
\begin{equation}
    p_{|S}\sim\frac{\Gamma(\alpha X+N)}{V^{\alpha X+N-1}}\EWD{u}{u\in\R_+^S\land\EWD{\sum i}{i\in S}{u(i)}=V}{\EWD{\prod i}{i\in S}{\frac{u(i)^{a(i)+c(i)-1}}{\Gamma(a(i)+c(i))}}}\;.
\end{equation}
If $Z$ is drawn first (after drawing $c$) then the densities for $p_{|S}$ and $W$ become
\begin{equation}
    p_{|S},W\sim\frac{\Gamma(\alpha+N)}{Z^{\alpha+N-1}}\EWD{u,t}{u\in\R_+^S\land t\in\R_+\land\EWD{\sum i}{i\in S}{u(i)}+t=Z}{\EWD{\prod i}{i\in S}{\frac{u(i)^{a(i)+c(i)-1}}{\Gamma(a(i)+c(i))}}\frac{W^{\alpha Y-1}}{\Gamma(\alpha Y)}}\;.
\end{equation}

\section*{Appendix: Approximating the Rao-Blackwellized estimator}

The Rao-Blackwell procedure provides coefficients $v$ concentrated on $S$ for which
\begin{equation}
\begin{split}
    \bar f&=\frac1{Z}\EWD{\sum i}{i\in D}{p(i)f(i)}\\
    &\simeq\frac1{N}\EWD{\sum i}{i\in S}{v(i;S)f(i)}\;.
\end{split}
\end{equation}
Inverse probability weighting provides a coefficient vector $\pi$, depending on $Z$, for which
\begin{equation}
    \bar f\simeq\frac1{Z}\EWD{\sum i}{i\in S}{\frac{p(i)}{\pi(i;Z)}f(i)}\;.
\end{equation}
Both estimates are unbiased (if the correct value of $Z$ is used in the inverse probability weighting).
Therefore
\begin{equation}
    \frac1{N}v(i;S)\simeq\frac1{Z}\frac{p(i)}{\pi(i;Z)}\qquad\textrm{($i\in S$)}\;.
\end{equation}
This looks odd as the left hand side depends explicitly on $S$ but not on $Z$ whereas $p/\pi$ on the right hand side depends on $Z$ but has no explicit dependence on $S$.
Of course, it is only an identity in expectation.

Recall
\begin{equation}
\begin{split}
    v(i;S)&=p(i)\frac\partial{\partial p(i)}\log(F_N(S,p_{|S}))\qquad\textrm{($i\in D$)}\\
    F_N(S,p_{|S})&=\EWD{\sum m_{|S}}{[1\leq m_{|S}]\land\EWD{\sum i}{i\in S}{m(i)}=N}{N!\EWD{\prod i}{i\in S}{\frac{p(i)^{m(i)}}{m(i)!}}}\;.
\end{split}
\end{equation}
Note that $[1\leq v_{|S}]$ and $[v_{|D\setminus S}=0]$.
We will develop the simplest asymptotic approximation and show its connection to the variable-$N$ (Poisson sampling) setting.
Introduce the exponential generating function
\begin{equation}
\begin{split}
    G_\lambda(S,p_{|S})&=\EWD{\sum N}{0\leq N}{(\lambda^N/N!)F_N(S,p_{|S})}\\
    &=\EWD{\prod i}{i\in S}{\EWD{\sum k}{1\leq k}{(\lambda p(i))^k/k!}}\\
    &=\left\langle\prod i:i\in S:e^{\lambda p(i)}-1\right\rangle\;.
\end{split}
\end{equation}
Conversely, using $\Coef$ as the ``coefficient of'' symbolic operator,
\begin{equation}
    F_N(S,p_{|S})=N!\Coef[\lambda^N]G_\lambda(S,p_{|S})\;.
\end{equation}
Ref.~\cite{FS2009}, Ch.~8, describes the saddle point method for development of asymptotic expansions for such expressions.
Here I am only concerned with the leading term, and only for the logarithmic gradient of $F_N$ and not its magnitude.
Therefore $\lambda$ is determined to solve $(\partial/\partial\lambda)\log G_\lambda(S,p_{|S})=N/\lambda$, or
\begin{equation}
    N=\EWD{\sum i}{i\in S}{\frac{\lambda p(i)}{1-e^{-\lambda p(i)}}}\;.
\end{equation}
(If $N=M$ then one must let $\lambda\to0$.)
The weights $v(i;S)$ ($i\in S$) are then approximated by
\begin{equation}
\begin{split}
    v(i;S)&\simeq p(i)\frac\partial{\partial p(i)}\log(G_\lambda(S,p_{|S}))\\
    &=\frac{\lambda p(i)}{1-e^{-\lambda p(i)}}\;.
\end{split}
\end{equation}
The same result is obtained for inverse probability weighting with Poisson sampling.
In that case one needs the distribution of $\EWD{i}{i\in D}{c(i)}$ given $S$; thus
\begin{equation}
\begin{split}
    c(i)&\sim\EWD{k}{1\leq k}{(\lambda p(i))^k/k!}/(e^{\lambda p(i)}-1)\qquad\textrm{($i\in S$)}\\
    c(i)&=0\qquad\textrm{($i\in D\setminus S$)}
\end{split}
\end{equation}
Conditioned on $i\in S$ we have $\Expec(c(i))=\lambda p(i)/(1-e^{-\lambda p(i)})$
The resulting Rao-Blackwellized estimator for $\bar f$ is identical to what was obtained via the exponential generating function.
\begin{equation}
    \bar f\simeq\frac1{N}\EWD{\sum i}{i\in S}{\frac{\lambda p(i)f(i)}{1-e^{-\lambda p(i)}}}\;.
\end{equation}
The derivation of this approximation via Poisson sampling invites the identification $\lambda=N/Z$.

\section*{Appendix: Expectation values}

As the underlying probability distribution $L0$ is in pure product form, relevant expectation values (before sampling) follow directly from their point values.
\begin{equation}
\begin{split}
    \Expec(M)&=\EWD{\sum i}{i\in D}{\Probab(1\leq c(i))}\\
    \Expec(N)&=\EWD{\sum i}{i\in D}{\Expec(c(i))}\\
    \Expec(U)&=\EWD{\sum i}{i\in D}{x(i)\Expec(\log p(i)\Char(1\leq c(i)))}\\
    \Expec(V)&=\EWD{\sum i}{i\in D}{\Expec(p(i)\Char(1\leq c(i)))}\\
    \Expec(W)&=\EWD{\sum i}{i\in D}{\Expec(p(i)\Char(c(i)=0))}\\
    \Expec(X)&=\EWD{\sum i}{i\in D}{x(i)\Probab(1\leq c(i))}\\
    \Expec(Y)&=\EWD{\sum i}{i\in D}{x(i)\Probab(c(i)=0)}\\
    \Expec(Z)&=\EWD{\sum i}{i\in D}{\Expec(p(i))}\;.
\end{split}
\end{equation}
These expressions depend on the vector $x$ and the parameters $(\alpha,b,\lambda)$, but we suppress that in the notation.
The involved pointwise expectations include:
\begin{equation}
\begin{split}
    \Probab(c(i)=0)&=(1+\lambda/b)^{-a(i)}\\
    \Probab(1\leq c(i))&=1-(1+\lambda/b)^{-a(i)}\\
    \Expec(c(i))&=\lambda a(i)/b\;.
\end{split}
\end{equation}
The right hand sides are consistent with $\Probab(1\leq c(i))\le\Expec(c(i))$.
Furthermore
\begin{equation}
\begin{split}
    \Expec(\log p(i)\Char(1\leq c(i)))&=\frac{b^{a(i)}}{\Gamma(a(i))}\EWD{\sumint t,k}{0\leq t\land1\leq k}{(\log t)t^{a(i)-1}e^{-(b+\lambda)t}\frac{(\lambda t)^k}{k!}\,\d t}\\
    &=\frac{b^{a(i)}}{\Gamma(a(i))}\frac\partial{\partial a(i)}(\Gamma(a(i))(b^{-a(i)}-(b+\lambda)^{-a(i)}))\\
    &=(\psi(a(i))-\log b)(1-(1+\lambda/b)^{-a(i)})+\log(1+\lambda/b)(1+\lambda/b)^{-a(i)}\;.
\end{split}
\end{equation}
And
\begin{equation}
\begin{split}
    \Expec(p(i)\Char(1\leq c(i)))&=\frac{\Gamma(a(i)+1)(b^{-a(i)-1}-(b+\lambda)^{-a(i)-1})}{\Gamma(a(i))b^{-a(i)}}\\
    &=\frac{a(i)}{b}(1-(1+\lambda/b)^{-a(i)-1})\\
    \Expec(p(i)\Char(c(i)=0))&=\frac{\Gamma(a(i)+1)(b+\lambda)^{-a(i)-1}}{\Gamma(a(i))b^{-a(i)}}\\
    &=\frac{a(i)}{b}(1+\lambda/b)^{-a(i)-1}\\
    \Expec(p(i))&=\frac{a(i)}{b}\;.
\end{split}
\end{equation}

Consider next the corresponding expectation values when $S\subseteq D$ is given, but without observing $c_{|S}$ or $p_{|S}$; it is only known that $i\in S\Leftrightarrow1\le c(i)$.
In this case $M$, $X$ and $Y$ are known.
Further:
\begin{equation}
\begin{split}
    \Expec(N|S)&=\EWD{\sum i}{i\in S}{\Expec(c(i)|1\leq c(i))}\\
    \Expec(U|S)&=\EWD{\sum i}{i\in S}{x(i)\Expec(\log p(i)|1\leq c(i))}\\
    \Expec(V|S)&=\EWD{\sum i}{i\in S}{\Expec(p(i)|1\leq c(i))}\\
    \Expec(W|S)&=\EWD{\sum i}{i\in D\setminus S}{\Expec(p(i)|c(i)=0)}\\
    \Expec(Z|S)&=\Expec(V|S)+\Expec(W|S)\;.
\end{split}
\end{equation}
As before we suppress the dependence on the parameters $\alpha$, $b$ and $\lambda$.
To obtain the pointwise expectations conditioned on $i\in S$ or on $i\in D\setminus S$ one divides by the appropriate probability.
Thus
\begin{equation}
\begin{split}
    \Expec(c(i)|1\leq c(i))&=\frac{\lambda a(i)/b}{1-(1+\lambda/b)^{-a(i)}}\\
    \Expec(\log p(i)|1\leq c(i))&=\psi(a(i))-\log b+\frac{\log(1+\lambda/b)}{(1+\lambda/b)^{a(i)}-1}\\
    \Expec(p(i)|1\leq c(i))&=\frac{a(i)}{b}\frac{1-(1+\lambda/b)^{-a(i)-1}}{1-(1+\lambda/b)^{-a(i)}}\\
    \Expec(p(i)|c(i)=0)&=\frac{a(i)}{b+\lambda}\;.    
\end{split}
\end{equation}

Finally, for the observable $N$ it is also of interest to consider the expectation after observing both $S$ and $p_{|S}$.
\begin{equation}
    \Expec(N|S,p_{|S})=\EWD{\sum i}{i\in S}{\Expec(c(i)|1\leq c(i),p(i))}\;.
\end{equation}
The pointwise expectation is
\begin{equation}
    \Expec(c(i)|1\leq c(i),p(i))=\frac{\lambda p(i)}{1-e^{-\lambda p(i)}}\;.
\end{equation}

\section*{Appendix: Supporting the Bayesian and the Profile Likelihood analysis}

The following analysis relies on positivity and monotonicity properties of the trigamma function and some related functions that follow from their representation as the Laplace transform of a positive function; i.e., their representation as a totally positive function.
Recall that for a totally positive function $f$ derivatives are of alternating sign: $[0\leq(-1)^kf^{(k)}]$ for $0\leq k$.
The representations are (for positive real $z$ for our purposes)
\begin{equation}
\begin{split}
    \psi'(z)&=\EWD{\sumint t}{0\leq t}{\frac{te^{-tz}}{1-e^{-t}}\,\d t}\\
    z\psi'(z)-1&=\EWD{\sumint t}{0\leq t}{\frac{e^t-1-t}{(e^t-1)(1-e^{-t})}e^{-tz}\d t}\\
    z^{-1}+1-z\psi'(z)&=\EWD{\sumint t}{0\leq t}{\frac{t+e^{-t}-1}{(e^t-1)(1-e^{-t})}e^{-tz}\d t}\;.
\end{split}
\end{equation}
Function $f:\R_+\to\R$ is called subadditive if $\EWD{\forall x,y}{0\leq x\land0\leq y}{f(x+y)\leq f(x)+f(y)}$ and one way to show subadditivity of $f$ is to show that $\EWD{z}{0\leq z}{f(z)/z}$ is descending.
Thus the functions $\EWD{z}{0\leq z}{z^2\psi'(z)}$ and $\EWD{z}{0\leq z}{1-z^2\psi'(z)}$ are subadditive.

Bayesian estimation of the missing mass or partition function involves integration of $L4$ and of $L5$ over $\alpha$ with prior $1/\alpha$, and we must show that this integral is well defined.
Observe
\begin{equation}
\begin{split}
    \frac\partial{\partial\alpha}\log L4(..,W;\alpha)&=\psi(\alpha)-\EWD{\sum i}{i\in S}{x(i)\psi(\alpha x(i))}-Y\psi(\alpha Y)+U+Y\log W-\log(V+W)\\
    \frac{\partial^2}{\partial\alpha^2}\log L4(..,W;\alpha)&=\psi'(\alpha)-\EWD{\sum i}{i\in S}{x(i)^2\psi'(\alpha x(i))}-Y^2\psi'(\alpha Y)\;.
\end{split}
\end{equation}
By subadditivity, $\alpha^2(\partial^2/\partial\alpha^2)\log L4(..,W;\alpha)\leq0$ everywhere, so that $\EWD{\alpha}{0\leq\alpha}{L4(..,W;\alpha)}$ is logarithmically concave.
The limiting behavior of $(\partial/\partial\alpha)\log L4$ is
\begin{equation}
\begin{split}
    \frac\partial{\partial\alpha}\log L4(..,W;\alpha)&=\frac{M}\alpha+O(1)\qquad\hbox{($\alpha\to0$)}\\
    &=-\Delta+\frac{M}{2\alpha}+O(\alpha^{-2})\qquad\hbox{($\alpha\to\infty$)}\;.
\end{split}
\end{equation}
Putting this together it is seen that $\alpha^{-1}L4(..,W;\alpha)$ is integrable except in the singular case $\Delta=0$.

One finds similar expressions for $L5$
\begin{equation}
\begin{split}
    \frac\partial{\partial\alpha}\log L5(..;\alpha)&=\psi(\alpha)-\EWD{\sum i}{i\in S}{x(i)\psi(\alpha x(i))}+U-X\log V+X\psi(\alpha X+N)-\psi(\alpha+N)\\
    \frac{\partial^2}{\partial\alpha^2}\log L5(..;\alpha)&=\psi'(\alpha)-\EWD{\sum i}{i\in S}{x(i)^2\psi'(\alpha x(i))}+X^2\psi'(\alpha X+N)-\psi'(\alpha+N)\;.
\end{split}
\end{equation}
To show that $\EWD{\alpha}{0\leq\alpha}{L5(..;\alpha)}$ is logarithmically concave use subadditivity again to establish
\begin{equation}
    \frac{\partial^2}{\partial\alpha^2}\log L5(..;\alpha)\leq\psi'(\alpha)-X^2\psi'(\alpha X)+X^2\psi'(\alpha X+N)-\psi'(\alpha+N)\;.
\end{equation}
Then use the trigamma function identity
\begin{equation}
    \psi'(z+N)-\psi'(z)=\EWD{\sum k}{0\leq k<N}{(z+k)^{-2}}
\end{equation}
to establish $(\partial^2/\partial\alpha^2)\log L5(..;\alpha)\leq0$.
The limiting behavior of $(\partial/\partial\alpha)\log L5$ is
\begin{equation}
\begin{split}
    \frac\partial{\partial\alpha}\log L5(..;\alpha)&=\frac{M-1}\alpha+O(1)\qquad\hbox{($\alpha\to0$)}\\
    &=-\Delta_S+\frac{M-1}{2\alpha}+O(\alpha^{-2})\qquad\hbox{($\alpha\to\infty$)}\;.
\end{split}
\end{equation}
Putting this together it is seen that $\alpha^{-1}L5(..;\alpha)$ is integrable except in the singular case $\Delta_S=0$.

Profile likelihood estimation of the missing mass or partition function involves maximization of $L8$ over $\alpha$, and we need to study the function to make sure that it is a well-defined problem.
The subsequent mixed treatment requires maximization of $L9$, and it needs a similar analysis.
We calculate
\begin{equation}
\begin{split}
    \frac\partial{\partial\alpha}\log L8(..,W;\alpha)&=\log\alpha-\EWD{\sum i}{i\in S}{x(i)\psi(\alpha x(i))}-Y\psi(\alpha Y)+U+Y\log W-\log(V+W)\\
    \frac{\partial^2}{\partial\alpha^2}\log L8(..,W;\alpha)&=\alpha^{-1}-\EWD{\sum i}{i\in S}{x(i)^2\psi'(\alpha x(i))}-Y^2\psi'(\alpha Y)
\end{split}
\end{equation}
As when dealing with $L4$ and $L5$ for Bayesian estimation we establish logarithmic concavity.
\begin{equation}
\begin{split}
    \alpha^2\frac{\partial^2}{\partial\alpha^2}\log L8(..;\alpha)&=\EWD{\sum i}{i\in S}{\alpha x(i)(1-\alpha x(i)\psi'(\alpha x(i)))}+\alpha Y(1-\alpha Y\psi'(\alpha Y))\\
    &<0\;.
\end{split}
\end{equation}
The asymptotic behavior is
\begin{equation}
\begin{split}
    \frac\partial{\partial\alpha}\log L8(..;\alpha)&=\frac{M+1}\alpha+\log\alpha+O(1)\qquad\hbox{($\alpha\to0$)}\\
    &=-\Delta+\frac{M+1}{2\alpha}+O(\alpha^{-2})\qquad\hbox{($\alpha\to\infty$)}
\end{split}
\end{equation}
Putting this together it is seen that $\EWD{\alpha}{0\leq\alpha}{L8(..;\alpha)}$ has a unique maximum at finite argument except in the singular case $\Delta=0$ when one must let $\alpha\to\infty$.

Proceeding with similar analysis for $L9$
\begin{equation}
\begin{split}
    \frac\partial{\partial\alpha}\log L9(..;\alpha)&=\log\alpha-\EWD{\sum i}{i\in S}{x(i)\psi(\alpha x(i))}+U-X\log V+X\psi(\alpha X+N)-\psi(\alpha+N)\\
    \frac{\partial^2}{\partial\alpha^2}\log L9(..;\alpha)&=\alpha^{-1}-\EWD{\sum i}{i\in S}{x(i)^2\psi'(\alpha x(i))}+X^2\psi'(\alpha X+N)-\psi'(\alpha+N)\;.
\end{split}
\end{equation}
Logarithmic concavity:
\begin{equation}
\begin{split}
    \alpha^2\frac{\partial^2}{\partial\alpha^2}\log L9(..;\alpha)&\leq\alpha-(\alpha X)^2\psi'(\alpha X)+(\alpha X)^2\psi'(\alpha X+N)-\alpha^2\psi'(\alpha+N)\\
    &=\alpha-\alpha^2\psi'(\alpha)+\EWD{\sum k}{0\leq k<N}{(\alpha X)^2/(\alpha X+k)^2-\alpha^2/(\alpha +k)^2}\\
    &\leq\alpha(1-\alpha\psi'(\alpha))\\
    &<0\;.
\end{split}
\end{equation}
Asymptotic behavior:
\begin{equation}
\begin{split}
    \frac\partial{\partial\alpha}\log L9(..;\alpha)&=\frac{M}\alpha+\log\alpha+O(1)\qquad\hbox{($\alpha\to0$)}\\
    &=-\Delta_S+\frac{M}{2\alpha}+O(\alpha^{-2})\qquad\hbox{($\alpha\to\infty$)}\;.
\end{split}
\end{equation}
Therefore $\EWD{\alpha}{0\leq\alpha}{L9(..;\alpha)}$ has a unique maximum at finite argument except in the singular case $\Delta_S=0$ when one must let $\alpha\to\infty$.

\section*{Appendix: Maximum Likelihood Estimation of the missing mass}

The starting point for maximum likelihood estimation of the missing mass and its distribution is the likelihood $L3(S,p_{|S},c_{|S};\alpha,b,\lambda)$.
The strategy is to optimize $L3$ over $(\alpha,b,\lambda)$ given the observations.
Then the expected value of the missing mass is $\Expec(W)=\alpha Y/(b+\lambda)$ and its probability distribution is a Gamma distribution.
\begin{equation}
    W\sim\frac{(b+\lambda)^{\alpha Y}}{\Gamma(\alpha Y)}\EWD{t}{0\leq t}{t^{\alpha Y-1}e^{-(b+\lambda)t}}\;.
\end{equation}
The conditions for stationarity of $L3$ follow from
\begin{equation}
\begin{split}
    \frac\partial{\partial\alpha}\log L3(S,p_{|S},c_{|S};\alpha,b,\lambda)&=-\EWD{\sum i}{i\in S}{x(i)\psi(\alpha x(i))}+X\log b+U-Y\log(1+\lambda/b)\\
    \frac\partial{\partial b}\log L3(S,p_{|S},c_{|S};\alpha,b,\lambda)&=\frac\alpha{b}-V-\frac{\alpha Y}{b+\lambda}\\
    \frac\partial{\partial\lambda}\log L3(S,p_{|S},c_{|S};\alpha,b,\lambda)&=\frac{N}\lambda-V-\frac{\alpha Y}{b+\lambda}
\end{split}
\end{equation}
The conditions for stationarity wrt $b$ and $\lambda$ are solved by simple algebra.
One obtains $N=\lambda\alpha/b$ and then
\begin{equation}
\begin{split}
    b(\alpha)&=\frac\alpha{V}\frac{\alpha X+N}{\alpha+N}\\
    \lambda(\alpha)&=\frac{N}{V}\frac{\alpha X+N}{\alpha+N}\;.
\end{split}
\end{equation}
Let $L11(S,p_{|S},c_{|S};\alpha)=L3(S,p_{|S},c_{|S};\alpha,b(\alpha),\lambda(\alpha))$, then
\begin{equation}
\begin{split}
    &\frac\partial{\partial\alpha}\log L11(S,p_{|S},c_{|S};\alpha)\\
    &\quad=-\EWD{\sum i}{i\in S}{x(i)\psi(\alpha x(i))}-X\log((1+N/\alpha)V)+X\log(\alpha X+N)+U-Y\log(1+N/\alpha)\\
    &\quad=-\EWD{\sum i}{i\in S}{x(i)\psi(\alpha x(i))}-X\log V+X\log(\alpha X+N)+U-\log(1+N/\alpha)\\
    &\frac{\partial^2}{\partial\alpha^2}\log L11(S,p_{|S},c_{|S};\alpha)\\
    &\quad=-\EWD{\sum i}{i\in S}{x(i)^2\psi'(\alpha x(i))}+X^2/(\alpha X+N)+(N/\alpha)/(\alpha+N)\;.
\end{split}
\end{equation}
Now we would like to show that $L11$ is logarithmically concave (as we did for $L4$ and $L8$ earlier), but this is not going to fly.
Consider the limit $X\to0$ (and therefore $x(i)\to0$ for $i\in S$) of $(\partial^2/\partial\alpha^2)\log L11(..;\alpha)$:
\begin{equation}
    \frac{\partial^2}{\partial\alpha^2}\log L11(S,p_{|S},c_{|S};\alpha)\to-M/\alpha^2+(N/\alpha)/(\alpha+N)\qquad\textrm{($X\to0$)}\;.
\end{equation}
This might be expected to be negative (if $X\to0$ then expect $M\simeq N$), but it is not guaranteed to be negative.
Therefore I expect maximum likelihood estimation to be messy, and I discount it in favour of the profile likelihood approach.

\section*{Appendix: Moment matching strategies for determining the parameters $\alpha$, $b$, $\lambda$}

The starting point for moment matching is the likelihood $L3(S,p_{|S},c_{|S};\alpha,b,\lambda)$.
It looks natural to match the moments $(N,U,V)$ that are sufficient given $S$, but that still leaves plenty of freedom.
Having determined $(\alpha,b,\lambda)$ the expected value of the missing mass is $\Expec(W)=\alpha Y/(b+\lambda)$ and its probability distribution is a Gamma distribution as for maximum likelihood estimation.

\subsection*{Moment Matching A}

Determine parameters $\alpha$, $b$, $\lambda$ so that the observed $N$, $U$, $V$ match the priors $\Expec(N)$, $\Expec(U)$, $\Expec(V)$ before any sampling.
Then the parameters $\alpha$, $b$, $\lambda$ are determined to satisfy:
\begin{equation}
\begin{split}
    N&=\lambda\alpha/b\\
    U&=\EWD{\sum i}{i\in D}{x(i)\left((\psi(a(i))-\log b)(1-(1+\lambda/b)^{-a(i)})+\log(1+\lambda/b)(1+\lambda/b)^{-a(i)}\right)}\\
    &=-\log(b)+\log(b+\lambda)\EWD{\sum i}{i\in D}{x(i)(1+\lambda/b)^{-a(i)}}+\EWD{i}{i\in D}{x(i)\psi(a(i))(1-(1+\lambda/b)^{-a(i)})}\\
    V&=\EWD{\sum i}{i\in D}{\frac{a(i)}{b}\left(1-(1+\lambda/b)^{-a(i)-1}\right)}\\
    &=\frac\alpha{b}\left(1-\EWD{\sum i}{i\in D}{x(i)(1+\lambda/b)^{-a(i)-1}}\right)\;.
\end{split}
\end{equation}
(All the time, $a(i)=\alpha x(i)$.)
The equation for $N$ tells us $\lambda/b=N/\alpha$ and then the equation for $V$ can be solved to obtain
\begin{equation}
    b=\frac\alpha{V}\left(1-\EWD{\sum i}{i\in D}{x(i)(1+N/\alpha)^{-a(i)-1}}\right)\;.
\end{equation}
After that the equation for $U$ must determine $\alpha$, but it does not look appealing.

\subsection*{Moment Matching B}

Determine parameters $\alpha$, $b$, $\lambda$ so that the observed $N$, $U$, $V$ match the priors $\Expec(N|S)$, $\Expec(U|S)$, $\Expec(V|S)$ after $S$ has been determined (but without observing $p_{|S}$).
Then the parameters $\alpha$, $b$, $\lambda$ are determined to satisfy:
\begin{equation}
\begin{split}
    N&=\frac{\lambda\alpha}{b}\EWD{\sum i}{i\in S}{\frac{x(i)}{1-(1+\lambda/b)^{-a(i)}}}\\
    U&=\EWD{\sum i}{i\in S}{x(i)\left(\psi(a(i))-\log b+\frac{\log(1+\lambda/b)}{(1+\lambda/b)^{a(i)}-1}\right)}\\
    V&=\EWD{\sum i}{i\in S}{\frac{a(i)}{b}\frac{1-(1+\lambda/b)^{-a(i)-1}}{1-(1+\lambda/b)^{-a(i)}}}\;.
\end{split}
\end{equation}

\subsection*{Moment Matching C}

Determine the parameter $\lambda$ so that the observed $N$ matches $\Expec(N|S,p_{|S})$ and determine parameters $\alpha$, $b$ to match $\Expec(U|S)$ and $\Expec(V|S)$ as in strategy B.
Then the parameter $\lambda$ is determined to satisfy
\begin{equation}
    N=\EWD{\sum i}{i\in S}{\frac{\lambda p(i)}{1-e^{-\lambda p(i)}}}\;.
\end{equation}
For the pair $(\alpha,b)$ the previous equations for $U$ and $V$ are used.

If any of these moment matching strategies should be pursued then it would be Strategy C, I think; it makes the most use of the available information by conditioning $\Expec(N)$ on both $S$ and $p_{|S}$ and conditioning $\Expec(U)$ and $\Expec(V)$ on $S$.
But in the end I am not getting any inspiration from these equations.

\end{document}